\documentclass[a4paper,11pt]{article}

\usepackage{url}
\usepackage{amssymb,amsmath,amsthm,amsfonts}

\newtheorem{theoreme}{Theorem}[section]

\newtheorem{prop}[theoreme]{Proposition}
\newtheorem{lemma}[theoreme]{Lemma}

\def\ssum{\sum\limits}

\font\tenCal=cmsy10
\newfam\Calfam
\textfont\Calfam=\tenCal

\def\CC{{\mathbb C}}
\newenvironment{demo}{\begin{trivlist}\item[]{\bf
Proof.}}{\hfill$\square$\end{trivlist}}
\def\dim{\mathop{\hbox{\textrm dim}}\nolimits}
\def\Ker{\mathop{\hbox{\textrm Ker}}\nolimits}

\def\Hoch{H\!H}
\def\Pois{H\!P}
\def\ssum{\sum\limits}
\def\blanc{\ \ \ \ \ \ \ \ \ \ }
\def\esp{\ \ \ \ \ }

\def\Ker{\mathop{\hbox{\textrm Ker}}\nolimits}

\def\dim{\mathop{\hbox{\textrm dim}}\nolimits}

\def\CC{{\mathbb C}}
\def\blanc{\ \ \ \ \ \ \ \ \ \ }
\def\esp{\ \ \ \ \ }

\title{Poisson (co)homology of truncated polynomial algebras in two variables}

\author{St\'ephane Launois\thanks{This research was supported by a Marie Curie European Reintegration Grant within the $7^{\mbox{th}}$ European Community Framework
  Programme.} $ $ and Lionel Richard\thanks{Research supported by  EPSRC Grant EP/D034167/1. L.R. wishes to thank the IMSAS of the University of Kent (UK) for its kind hospitality during his stays in Canterbury in September 2007 and April 2008.}}

\date{}

\begin{document}

\maketitle
\begin{abstract}
We study the Poisson (co)homology of the algebra of truncated polynomials in two variables viewed as the semi-classical limit of a quantum complete intersection studied by Bergh and Erdmann. We show in particular that the Poisson cohomology ring of such a Poisson algebra is isomorphic to the Hochschild cohomology ring of the corresponding quantum complete intersection.
\end{abstract}


\section{Introduction}
\label{}
Given a Poisson algebra, its Poisson cohomology--as introduced by Lichnerowicz in \cite{Lichn}--provides important informations about the Poisson structure (the
Casimir elements are reflected by the degree zero cohomology,
Poisson derivations modulo Hamiltonian derivations by the degree
one,...). It also plays a crucial role in the study of deformations of the Poisson structure.
A classical problem in algebraic deformation is to compare the Poisson (co)homology of a Poisson algebra with the Hochschild (co)homology of its deformation.
Although these homologies are known to behave similarly in smooth cases (see e.g. \cite{brylinski,K,Kont,LMP}), the singular case seems more complicated to deal with.
There are many examples in which the trace groups already do not match in \cite{AF}. Even though the dimension of the homology spaces in degree zero match, as for Kleinian singularities in \cite{AL98}, it might not be the case in higher degree (see \cite{AFLS} and \cite{Pichereau}).

Our aim in this paper is to provide an example of a singular Poisson algebra
 such that its Poisson  cohomology ring is isomorphic as a graded commutative algebra to the Hochschild  cohomology ring of a natural deformation.
Namely for two integers $a,b\geq 2$ we consider the truncated polynomial algebra $\Lambda(a,b):=\CC[X,Y]/( X^a,Y^b)$ with the Poisson bracket given by  $\{X,Y\}=XY$.
This algebra is the semi-classical limit of the quantum complete intersection $\Lambda_{q}(a,b)$ studied by Bergh and Erdmann in \cite{BE}. The homological properties of this class of noncommutative finite-dimensional algebras have been extensively studied recently (see for instance \cite{BGMS,Bergh,BE,BO}).
In particular it is proved in \cite{BGMS} that the generic quantum exterior algebra (corresponding to the case $a=b=2$) provides a counter-example to Happel's question.
That is, the global dimension of $\Lambda_{q}(2,2)$ is infinite whereas its Hochschild cohomology groups vanish in  degree greater than 2.
Recently the Hochschild cohomology ring of $\Lambda_{q}(a,b)$ has been described in \cite{BE} when $q$ is not a root of unity. It is a five-dimensional  graded algebra isomorphic to the fibre product ${\mathcal F}:=\CC[U]/(U^2) \times_{\CC} \CC\langle V,W\rangle / (V^2, VW+WV, W^2)$, with $U$ in degree zero and $V$ and $W$ in degree one. 
In this paper we obtain the same description for the Poisson cohomology ring of $\Lambda(a,b)$, giving rise to

\medskip

\begin{theoreme} \label{thm:coh}
For all $a,b\geq 2$ one has
$\Pois^*(\Lambda(a,b)) \simeq  {\mathcal F} \simeq \Hoch^*(\Lambda_{q}(a,b))$ as graded commutative algebras, where $\Pois^*(\Lambda(a,b))$ is endowed with the Poisson cup product induced by the exterior product on skew-symmetric multiderivations.
\end{theoreme}

\medskip

Note that $\Lambda(a,b)$ is not unimodular (see for instance \cite{Kosmann} for a recent survey on the modular class of a Poisson manifold). Indeed, there is no Poincar\'e duality  between $\Pois^k(\Lambda(a,b))$ and $\Pois_{2-k}(\Lambda(a,b))$. For instance $\Pois^2(\Lambda(a,b))$ has dimension 1 whereas $\Pois_0(\Lambda(a,b))$ has dimension $a+b-1$.
 One cannot even expect a twisted Poincar\'e duality as in \cite{LMP} (see the paragraph below Lemma \ref{lem:trace}). However, the Nakayama automorphism $\nu$ coming from the Frobenius structure of the quantum complete intersection $\Lambda_q(a,b)$ (see \cite{BE}) allows us to construct a Poisson module $M_\nu$ such that 
 
 \medskip
 
 \begin{theoreme}\label{thm:dual}
$\Pois^k(\Lambda(a,b))\simeq \Pois_k(\Lambda(a,b),M_\nu)$ for all nonnegative integer  $k$.
\end{theoreme}

\section{Poisson cohomology of truncated polynomials}

The algebra $\Lambda(a,b)$ has dimension $ab$, with basis $\{e_{ij}:=X^iY^j\ |\ 0\leq i\leq a-1,\ 0\leq j\leq b-1\}$.
One easily checks that the following formulas hold in $\Lambda(a,b)$ for all $i,j$:
\begin{equation}\label{bracket}
\{X^iY^j,X\}=-jX^{i+1}Y^j; \ \ \ \ \{X^iY^j,Y\}=iX^{i}Y^{j+1}.
\end{equation}

Recall that the complex computing the Poisson cohomology is $(\chi^k,\delta_k)$, with $\chi^k$ the space of skew-symmetric $k$-derivations of $\Lambda(a,b)$. Interestingly, these spaces vanish for $k\geq 3$, since $\Lambda(a,b)$ is 2-generated. 
The Poisson coboundary operator $\delta_k : \chi^k \rightarrow \chi^{k+1}$ is defined by
\begin{eqnarray*}
\delta_k(P)(f_0,\dots,f_k) & := &  \sum_{i=0}^k(-1)^i \left\{ f_i, P(f_0,\dots,\widehat{f_i},\dots,f_k)\right\}\\
& &  +\sum_{0 \leq i < j \leq k} (-1)^{i+j}P\left( \{f_i,f_j\},f_0, \dots,\widehat{f_i},\dots,\widehat{f_j}, \dots, f_k \right)
\end{eqnarray*}
for all $P \in \chi^k$ and $f_0,\dots,f_k \in \Lambda(a,b)$. It is easy to check that $\delta_k(P)$ belongs indeed to $\chi^{k+1}$ and that
$\delta_{k+1} \circ \delta_k =0$. The $k^{\rm{th}}$  Poisson cohomology space of $\Lambda(a,b)$,  denoted by $\Pois^k(\Lambda(a,b))$, is the $k^{\rm{th}}$ cohomology space of this complex. 
The Poisson cohomology ring is the space $\Pois^*(\Lambda(a,b)):=\oplus_{k=0}^{\infty}\Pois^k(\Lambda(a,b))$. Endowed with the cup product induced by the exterior product on $\chi^*:=\oplus_k \chi^k$ it becomes a graded commutative algebra.
We start by describing the spaces $\Pois^k(\Lambda(a,b))$ for $k=0,1,2$ as it is clear that $\Pois^k(\Lambda(a,b))=0$ for $k\geq 3$.
 
 Let $\lambda=\sum \lambda_{ij}e_{ij}\in \Lambda(a,b)$, and assume it satisfies $\{\lambda,X\}=\{\lambda,Y\}=0$. Then  equations (\ref{bracket}) lead successively to $\lambda_{ij}=0$ for all $i\neq a-1$ and all $j\neq 0$, and then to $\lambda_{ij}=0$ for all $i\neq 0$ and $j\neq b-1$. Hence
 
\medskip


\begin{prop}
The Poisson centre $\Pois^0(\Lambda(a,b))$ is  equal to $\CC\oplus\CC X^{a-1}Y^{b-1}=\CC e_{0,0}\oplus\CC e_{a-1,b-1}$.
\end{prop}

\medskip

 Any derivation $d\in\chi^1$ is uniquely determined by the values of $d(X)=\sum \lambda_{ij}e_{ij}$ and $d(Y)=\sum \lambda'_{ij}e_{ij}$. Moreover, $d$ must satisfy the relations $d(X^a)=d(Y^b)=0$, that is 
$X^{a-1}d(X)=Y^{b-1}d(Y)=0$, 
 since $d$ is a derivation. From that one easily deduces that $\lambda_{0j}=\lambda'_{i0}=0$ for all $i,j$.  Hence
the space $\chi^1={\rm Der}(A,A)$ has  basis $\{d_{ij}\}\cup\{d'_{ij}\}$, where: 
 \begin{enumerate}
 \item for $1\leq i< a$ and $0\leq j<b$, the derivation $d_{ij}$ is defined by $d_{ij}(X)=X^iY^j$ and $d_{ij}(Y)=0$;
 \item for $0\leq i< a$ and $1\leq j<b$, the derivation $d'_{ij}$ is defined by $d'_{ij}(X)=0$ and $d'_{ij}(Y)=X^iY^j$.
\end{enumerate}
In particular, ${\rm dim}(\chi^1)=b(a-1)+a(b-1)$. 
 Let $d=\sum_{i\neq 0}\alpha_{i,j}d_{i,j} + \sum_{j\neq 0}\beta_{i,j}d'_{i,j}\in\chi^1$. Then $d\in{\rm Ker}\delta_1$ if and only if it satisfies $d\left(\{X,Y\}\right)=\{d(X),Y\}+\{X,d(Y)\}$, that is:
 $$\sum_{j \neq 0}\beta_{i,j}X^{i+1}Y^{j}+\sum_{i \neq 0}\alpha_{i,j}X^{i}Y^{j+1}=\sum_{i \neq 0} i\alpha_{i,j}X^{i}Y^{j+1}+\sum_{j \neq 0} j\beta_{i,j}X^{i+1}Y^{j}.$$
 Identifying the coefficients in front of $X^iY^j$ leads to

 \begin{equation}\label{eqderpoisson}
 d\in{\rm Ker}\delta_1 \ \iff \ (1-j)\beta_{i-1,j}+(1-i)\alpha_{i,j-1}=0 \ \ \ \forall \ 1\leq i\leq a-1,\ 1\leq j\leq b-1.
\end{equation}


\begin{prop}
 $\Pois^1(\Lambda(a,b))=\CC d_{1,0}\oplus\CC d'_{0,1}$.
\end{prop}

\smallskip

\begin{demo}
 Let $d=\sum_{i\neq 0}\alpha_{i,j}d_{i,j} + \sum_{j\neq 0}\beta_{i,j}d'_{i,j}\in{\rm Ker}\delta_1$. Set $\lambda:=\sum_{j \neq 0} \frac{\alpha_{i+1,j}}{j} X^iY^j\in \Lambda(a,b)$. 
From (\ref{bracket}) one deduces that $d_1=d+\{\lambda,-\}$ is a Poisson derivation 
satisfying $d_1(X)=\sum_{i\geq 1}\alpha_{i,0}X^i$. Then one deduces from (\ref{eqderpoisson}) 
that $\alpha_{i,0}=0$ for all $i\neq 1$, that is $d_1(X)=\alpha_{1,0}X$, 
and  $d_1(Y)=\sum_{i=0}^{a-2}\beta_{i,1}X^iY+\sum_{j= 1}^{b-1}\beta_{a-1,j}X^{a-1}Y^j$. 
Now set $\mu:=\sum_{i=1}^{a-2} \frac{\beta_{i,1}}{i} X^i + \sum_{j= 1}^{b-1} \frac{\beta_{a-1,j}}{a-1} X^{a-1}Y^{j-1} $ and $d_2:=d_1-\{\mu,-\}$.
From the construction we get that $d_2=\alpha_{1,0}d_{1,0}+\beta_{0,1}d'_{0,1}$, 
so that the images of $d_{1,0}$ and $d'_{0,1}$ span $\Pois^1(\Lambda(a,b))$. 
One deduces from (\ref{bracket}) that they actually form a basis of $\Pois^1(\Lambda(a,b))$.
\end{demo}

\medskip

The complex  computing the Poisson cohomology is vanishing after $\chi^2$, so we  use the Euler-Poincar\'e principle to compute the dimension of $\Pois^2(\Lambda(a,b))$. First 
 note that a skew-symmetric derivation $f\in\chi^2$ is determined by $f(X\wedge Y)=\sum c_{ij}e_{ij}$, with $aX^{a-1}f(X\wedge Y)=bY^{b-1}f(X\wedge Y)=0$, so that  $c_{0j}=c_{i0}=0$ for all $i,j$. 
Thus
 $\chi^2$ has dimension $(a-1)(b-1)$.
\medskip

\begin{prop}
 $\Pois^2(\Lambda(a,b))=\CC f_{1,1}$, with $f_{1,1}:X\wedge Y\mapsto XY$.
\end{prop}

\smallskip

\begin{demo}
 We first prove that $\Pois^2(\Lambda(a,b))$ has dimension 1. From the Euler-Poincar\'e principle we get 
$\dim (\Pois^2(\Lambda(a,b)))=\dim \chi^2-{\rm rg}\delta_1=(a-1)(b-1)-(\dim\chi^1-\dim\Ker\delta_1)=(a-1)(b-1)-\left[a(b-1)+b(a-1)-(\dim \Pois^1(\Lambda(a,b))+{\rm rg}\delta_0)\right]=1-ab+(2+\dim \Lambda(a,b)-\dim\Ker\delta_0)=3-ab+ab-2=1$. 
Now all that remains is to check that $X\wedge Y\mapsto XY$ is not a Poisson  coboundary. 
Any $P\in\chi^1$  satisfies $\delta_1(P)(X\wedge Y)=\{X,P(Y)\}-\{Y,P(X)\} -P(XY)$. 
Moreover one has $P(XY)=P(X)Y+P(Y)X$, and it results straight from formulas (\ref{bracket}) 
that one cannot have $\delta_1(P)(X\wedge Y)=XY$.
\end{demo}

\medskip

\begin{trivlist}\item[]{\bf
Proof of Theorem \ref{thm:coh}.}
We have already proved that the graded commutative algebra $\Pois^*(\Lambda(a,b))$ is five-dimensional with basis $(e_{0,0}, e_{a-1,b-1}, d_{1,0}, d'_{0,1}, f_{1,1})$ with degree respectively (0,0,1,1,2).
One can easily check that $e_{a-1,b-1}$ annihilates $d_{1,0}$, $d'_{0,1}$ and $f_{1,1}$, and that $d_{1,0}\smile d'_{0,1}=f_{1,1}$.
Thus we have $\Pois^*(\Lambda(a,b))\simeq {\mathcal F}$, and one concludes using \cite[Theorem 3.3]{BE}.
\hfill$\square$\end{trivlist}

\medskip

Note that the basis $(e_{0,0}, e_{a-1,b-1}, d_{1,0}, d'_{0,1}, f_{1,1})$ of $\Pois^*(\Lambda(a,b))$ corresponds to the basis $(1,y^{b-1}x^{a-1},g,h,hg)$ considered in the proof of \cite[Theorem 3.3]{BE}.

\section{(Twisted) Poisson homology}

Let $M$ be a right Poisson module over $\Lambda(a,b)$ (the reader is referred to \cite[Section 3.1]{LMP} for the definition of a Poisson module).
We denote by $\{-,-\}_M:M\otimes_{\mathbb C} \Lambda(a,b)\longrightarrow M$ its external bracket.
 Then one defines a chain complex on the $\Lambda(a,b)$-modules $M\otimes_{\Lambda(a,b)} \Omega^k$, where $\Omega^k$ denotes the so-called K\"ahler differential  $k$-forms of $\Lambda(a,b)$, as follows.
The boundary operator $\partial_k:  M\otimes_{\Lambda(a,b)} \Omega^k \rightarrow  M\otimes_{\Lambda(a,b)} \Omega^{k-1}$ is defined by
$$\begin{array}{l}
\hspace{-.8cm} \partial_k (m \otimes da_1 \wedge \dots \wedge da_k )  =  \ssum_{i=1}^k (-1)^{i+1} \{m ,a_i\}_M \otimes da_1 \wedge \dots \wedge \widehat{da_i} \wedge \dots \wedge da_k \\
 \blanc \blanc \esp + \ssum_{1 \leq i < j \leq k} (-1)^{i+j} m \otimes d\{a_i,a_j\} \wedge da_1 \wedge \dots \wedge \widehat{da_i} \wedge \dots \wedge \widehat{da_j} \wedge \dots \wedge da_k ,
\end{array}$$
where we have removed the expressions under the hats in the previous sums and $d$ denotes the exterior differential. The homology of this
complex is denoted by $\Pois_{*}(\Lambda(a,b),M)$ and called the Poisson homology of the Poisson algebra $\Lambda(a,b)$ with
values in the Poisson module $M$.
Before computing the Poisson homology spaces for a specific module, we describe the spaces $\Omega^k$. By definition $\Omega^k$ is a $\Lambda(a,b)$-module generated by the wedge products of length $k$ of the $1$-differential forms $dX, dY$. In particular  $\Omega^k=0$ for $k\geq 3$.
In the remaining cases   the torsion coming from the relations $X^a=Y^b=0$ leads to the following spaces: $$\Omega^0=\Lambda(a,b) , \ \ \Omega^1=\bigoplus_{\substack{0\leq i\leq a-2 \\ 0\leq j\leq b-1}}\CC X^i Y^jdX\oplus \bigoplus_{\substack{0\leq i\leq a-1 \\ 0\leq j\leq b-2}}\CC X^iY^jdY,$$ and $$\Omega^2=\bigoplus_{\substack{0\leq i\leq a-2 \\ 0\leq j\leq b-2}}\CC X^iY^j dX\wedge dY.$$ Their dimensions are respectively $ab$, $(a-1)b+b(a-1)$, and $(a-1)(b-1)$.
As the bracket $\{\lambda,\mu\}$ belongs to the ideal generated by $XY$ for any $\lambda, \mu \in \Lambda(a,b)$, one easily checks that:

\medskip

\begin{lemma}\label{lem:trace}
$\Pois_0(\Lambda(a,b))$ has dimension $a+b-1$.
\end{lemma}

\medskip

As $a+b-1\geq 3$ for any $a,b\geq 2$, and $\Pois^2(\Lambda(a,b))$ has dimension 1,
this lemma shows that there is no Poincar\'e duality  between $\Pois^k(\Lambda(a,b))$ and $\Pois_{2-k}(\Lambda(a,b))$.

We may ask now if there is a twisted duality similar to the one obtained in \cite{LMP}. 
The Poisson algebra $\Lambda(a,b)$ is the semi-classical limit of the quantum complete intersection $\Lambda_q(a,b)$ which is the $\CC$-algebra generated by $x,y$ with relations $xy=qyx$, $x^a=0$, $y^b=0$ (see \cite{BE}).
Any diagonal automorphism $\sigma$ of $\Lambda_q(a,b)$ defined by $\sigma(x)=q^{\alpha}x$ and $\sigma(y)=q^{\beta}y$ gives rise to a Poisson module $M_\sigma$ of $\Lambda(a,b)$ via a semi-classical limit process (see \cite[Section 3.1]{LMP} for details). As a vector space $M_\sigma$ is equal to $\Lambda(a,b)$, and the external Poisson bracket is defined by:
\begin{equation}\label{form:sc-mod}
 \{X^iY^j, X\}_{M_\sigma}=-(j+\alpha) X^{i+1}Y^j \ \ \ \mbox{ and }\ \ \ \{X^iY^j, Y\}_{M_\sigma} =(i-\beta) X^{i}Y^{j+1}.
\end{equation}
From these formulas one easily deduces that if $\alpha<-b+1$ and $\beta>a-1$, then $\Pois_0(\Lambda(a,b),M_\sigma)=\CC$.
So one might be tempted to use such a Poisson module to restore a twisted Poincar\'e duality 
between $\Pois^k(\Lambda(a,b))$ and $\Pois_{2-k}(\Lambda(a,b),M_\sigma)$ as in \cite[Theorem 3.4.2]{LMP}. Unfortunately, these hypotheses on $\alpha$ and $\beta$ also lead to $\Pois_2(\Lambda(a,b),M_\sigma)=0$, so that $\Pois_2(\Lambda(a,b),M_\sigma)$ is not isomorphic to $\Pois^0(\Lambda(a,b))$.

The Nakayama automorphism $\nu$ of $\Lambda_q(a,b)$ coming from the Frobenius algebra structure of $\Lambda_q(a,b)$ is defined by $\nu(x)=q^{1-b}x$ and $\nu(y)=q^{a-1}y$.
This automorphism was used in \cite[Section 3]{BE} to link the twisted Hochschild homology and the Hochschild cohomology of $\Lambda_q(a,b)$ in each degree. We end this paper by computing the dimensions of the Poisson homology spaces of $\Lambda(a,b)$ with value in the Poisson module $M_\nu$ corresponding to the Nakayama automorphism $\nu$.

\medskip

\begin{prop}\label{prop:naka}
The twisted Poisson homology spaces $\Pois_k(\Lambda(a,b),M_\nu)$ have dimension 2,2,1 for $k=0,1,2$ respectively,  and are null if $k\geq 3$.
\end{prop}

\smallskip

\begin{demo}
 From $\partial_1(X^iY^j\otimes dX)=\{X^iY^j, X\}_{M_\nu} =-(j-b+1) X^{i+1}Y^j$ and $\partial_1(X^iY^j\otimes dY)=\{X^iY^j, Y\}_{M_\nu} =(i-a+1) X^{i}Y^{j+1}$ one easily gets that $\Pois_0(\Lambda(a,b),M_\nu)=\Lambda(a,b)/{\rm Im}\partial_1$ is generated as a $\CC$-vector space by the classes of $1$ and $X^{a-1}Y^{b-1}$ modulo Im$\partial_1$.

 We compute $\partial_2(X^iY^j\otimes dX\wedge dY)=\{X^iY^j, X\}_{M_\nu}\otimes dY - \{X^iY^j, Y\}_{M_\nu}\otimes dX-X^iY^j\otimes(XdY+YdX)=-(j-b+2) X^{i+1}Y^j\otimes dY-(i-a+2) X^{i}Y^{j+1}\otimes dX$, from which one easily sees that ${\rm Ker}\partial_2=\Pois_2(\Lambda(a,b),M_\nu)$ is the $\CC$-vector space generated by $X^{a-2}Y^{b-2}\otimes dX\wedge dY$.

 We use the Euler-Poincar\'e principle to conclude. There is only one dimension missing, all the other ones have the same value as for the Poisson cohomology of $\Lambda(a,b)$, so we end up with the same value for dim$\Pois_1(\Lambda(a,b),M_\nu)$, that is 2.
\end{demo}

\medskip

We deduce from Proposition \ref{prop:naka} and Theorem \ref{thm:coh} a duality similar to the one appearing in \cite{BE}, that is $\Pois^k(\Lambda(a,b))\simeq \Pois_k(\Lambda(a,b),M_\nu)$ for all nonnegative integer  $k$ as claimed in Theorem \ref{thm:dual}.





\bigskip

\noindent
St\'ephane Launois\\
Institute of Mathematics, Statistics and Actuarial Science\\
 University of Kent\\
  Canterbury CT2 7NF, UK\\
   E-mail: {S.Launois@kent.ac.uk}

\bigskip

\noindent
Lionel Richard\\
School of Mathematics\\
 University of Edinburgh\\
  Edinburgh EH9 3JZ, UK\\
   E-mail: Lionel.Richard@ed.ac.uk.

\end{document}